\documentclass[12pt]{article}
\usepackage[cp1251]{inputenc}
\usepackage[russian]{babel}
\usepackage[T2B]{fontenc}
\usepackage[pdftex, colorlinks, urlcolor=blue, pdfborder={0 0 0 [0 0]}]{hyperref}
\usepackage{amsfonts,amsmath,amssymb,theorem}
\usepackage{array}

\newtheorem{theorem}{Теорема}
\newtheorem{lemma}{Лемма}
\newtheorem{cor}{Следствие}
\newtheorem{prp}{Утверждение}
\newtheorem{remark}{Замечание}
\newtheorem{dfn}{Определение}

\newenvironment{proof}
{\vspace{1pt}\par{\sl%
Proof.\,\ }}%
{\noindent\vspace{1pt}}

\title{Об одном из обобщений формулы Ито-Вентцеля для случая нецентрированной пуассоновской меры и стохастическом первом интеграле }

\author{Е. В. Карачанская\\
Тихоокеанский госуниверситет, Хабаровск}

\begin{document}
\maketitle


{формула Ито-Вентцеля,  нецентрированная пуассоновская мера, ядро интегрального инварианта, стохастический первый интеграл}

\begin{abstract}
Строится обобщение формулы Ито-Вентцеля для системы обобщенных СДУ Ито с нецентрированной мерой на основе стохастического ядра интегрального преобразования. Построена система обобщенных СДУ, решение которой -- ядро интегрального инварианта, связанного с решением обобщенного СДУ с нецентрированной мерой. Введено понятие стохастического первого интеграла системы обобщенных СДУ с нецентрированной мерой и определены условия, при выполнении которых случайная функция является первым интегралом заданной системы обобщенных СДУ.
\end{abstract}

\section*{Введение}
В {\rm{\cite{D_78,D_89,D_02}}} был предложен специфический подход
к стохастическому интегрированию, основанный на том, что
производится отображение не точки, а  ее окрестности по некоторому
традиционного, основанного на исчислении Маллявена. Переход к
стохастическим уравнениям для конкретных реализаций происходит на
основе сопоставления динамики элемента пространственной структуры,
маркированной числом, c динамикой реализации траектории, для
которой это число является начальным значением. Эта возможность
обусловлена близостью, в вероятностном смысле, траекторий,
стартовавших в области исходной окрестности на каждом из временных
сечений.

В настоящее время является актуальным исследование открытых систем,
которые характеризуются наличием в них случайных возмущений, как непрерывных, описываемых винеровскими процессами, так и скачкообразных, описываемых пуассоновскими процессами. При этом, в условиях случайных возмущений, некоторые характеристики системы должны сохраняться. Т. е. возникают задачи построения программных управлений в открытых системах. Для случая, когда возмущение вызвано только винеровскими процессами, имеется несколько видов алгоритмов построения программных управлений. Наличие пуассоновской составляющей в случайном возмущении до сих пор не позволяло строить программное управление такими системами. Построенное обобщение формулы Ито-Вентцеля эту проблему снимает \cite{11_KchUpr}.

\section{Стохастический объем}

Пусть ${\bf x}(t)$ -- динамический процесс, определенный на
$\mathbb{R}^{n}$, являющийся решением системы стохастических
дифференциальных уравнений
\begin{equation}\label{Ayd01.1}
\begin{array}{c}
  dx_{i}(t)= \displaystyle a_{i}\Bigl(t;{\bf x}(t) \Bigr)\, dt
  +\sum\limits_{k=1}^{m}
b_{i,k}\Bigl(t;{\bf x}(t)\Bigr)\, dw_{k}(t) +
\int\limits_{\mathbb{R}(\gamma)}g_{i}\Bigl(t;{\bf x}(t);\gamma\Bigr)\nu(dt;d\gamma), \\
  {\bf x}(t)={\bf x}\Bigl(t;{\bf x}(0) \Bigr) \Bigr|_{t=0}={\bf
  x}(0), \ \ \ \ \ \ i=\overline{1,n}, \ \ \ \ \ \ t\geq 0,
\end{array}
\end{equation}
где ${\bf w} (t)$ -- $m$-мерный винеровский процесс,
$\nu(t;\Delta\gamma)$ -- однородная по $t$ нецентрированная мера
Пуассона.

Относительно коэффициентов ${ a}(t;{\bf x})$, ${ b}(t;{\bf x})$  и
${g}(t;{\bf x};\gamma)$ будем предполагать, что они ограничены и
непрерывны вместе со своими производными
$
\displaystyle \frac{\partial a_{i}(t;{\bf x})}{\partial x_{j}}$,\,
$\displaystyle \frac{\partial b_{i,k}(t;{\bf x})}{\partial
x_{j}}$,\, $\displaystyle \frac{\partial g_{i}(t;{\bf
x};\gamma)}{\partial x_{j}}$,\, $i,j=\overline{1,n}
$
по совокупности переменных $(t;{\bf x};\gamma)$.

Этих условий достаточно для существования случайных величин
\begin{equation*}
\textit{J}_{i,j}(t)= \mathop {\operatorname{l}
.i.m.}\limits_{\left|\Delta \right|= \Delta_{j}\rightarrow
0}\displaystyle \frac{x_{i}\Bigl(t; {\bf x}(0)+\Delta
\Bigr)-x_{i}\Bigl(t; {\bf x}(0)\Bigr)}{\Delta_{j}},
\end{equation*}
которые могут быть найдены как решения системы (\ref{Ayd01.1}) и
системы стохастических  уравнений {\rm{\cite[c.\,287]{GS_68}}}
\begin{equation}\label{Ayd2.2}
\begin{array}{c}
  d\textit{J}_{i,j}(t)=\textit{J}_{l,j}(t)\, \biggl[\displaystyle
\frac{\partial a_{i}(t;{\bf x}(t))}{\partial x_{l}(0)}\, dt +
\beta \,\displaystyle \frac{\partial b_{i,k}(t;{\bf
x}(t))}{\partial
x_{l}(0)}\, dw_{k}(t) + \displaystyle \int\limits_{\mathbb{R}(\gamma)}
\frac{\partial g_{i}(t;{\bf x}(t))}{\partial x_{l}(0)}\nu(dt;d\gamma) \biggr], \\
\textit{J}_{i,j}(t)=\textit{J}_{i,j}\Bigl(t; {\bf x}(0) \Bigr)
\Bigr|_{t=0}=\delta_{i,j}=\left\{
\begin{array}{ll}
  0, & i\neq j,  \\
  1, & i=j.
\end{array}
  \right.
\end{array}
\end{equation}

Пусть $\textit{J}(t)= \textit{J}\Bigl(t; {\bf x}(0)\Bigr)=\det\Bigl(\textit{J}_{ij}\Bigl(t; {\bf x}(0)\Bigr)\Bigr)$, $i,j=\overline{1,n}$. Тогда случайную величину $\textit{J}(t)$ можно рассматривать как якобиан перехода от ${\bf x}(0)$ к ${\bf x}(t;{\bf x}(0))$.

При помощи случайных величин $\textit{J}_{i,j}(t)$ можно построить
всевозможные случайные многообразия
${\cal M}\Bigl(t; {\bf x}(t;\lambda ) \Bigr)$ ($\lambda\in {\cal
D} \subset \mathbb{R}^{r}, \  {\bf x}(t
;\lambda )\in\mathbb{R}^{n}, \  r \leq n $), сопоставив их с
конфигурационными многообразиями исходного заданного многообразия
${\cal M}\Bigl(0; {\bf x}(0;\lambda ) \Bigr)$ ($\lambda\in {\cal
D} \subset \mathbb{R}^{r}, \  {\bf x}(0;\lambda )={\bf
x}(\lambda )\in\mathbb{R}^{n}, \ r \leq n $).

\begin{dfn}\label{Aydo2}{\rm{\cite{D_89}}} Элементом стохастического объема,
порожденного случайным  процессом ${\bf x}\Bigl(t; {\bf x}(0;
\lambda)\Bigr)$, назовем структуру
$
d\Gamma(t)=\textit{J}(t) \,d\Gamma\Bigl({\bf x}(0)\Bigr)
$,
где $\textit{J}(t)$ -- якобиан, построенный из элементов
$\textit{J}_{ij}\Bigl(t; {\bf x}(0)\Bigr)$,
$
d\Gamma({\bf x})= \displaystyle \prod \limits_{i=1}^{n} dx_{i}$.
\end{dfn}

На основе определения (\ref{Aydo2})
в {\rm{\cite{D_89}}} вводится случайная структура
\begin{equation}\label{Ayd2.3}
\begin{array}{c}
  {\cal I}_{\Gamma(t)}(t)=\underbrace{\int \cdots
\int}\limits_{\Gamma (t)} f(t;{\bf z}) \, d \Gamma ({\bf z})
  =\underbrace{\int \cdots \int} \limits_{\Gamma (0)}f\Bigl(t;{\bf
x}(t;{\bf y})\Bigr) \, \textit{J}(t)\, d \Gamma ({\bf y}),
\end{array}
\end{equation}
где ${\bf x}(t;{\bf y}) $ -- случайный процесс со случайным начальным условием ${\bf y}$.
При этом интегрирование по случайному
объему $\Gamma(t)$ проводится на основе
интегральных сумм Римана, понятия среднеквадратического предела
для любых непрерывных по ${\bf z}$, в общем, случайных функций
$f(t;{\bf z})$. Структура \eqref{Ayd2.3} есть, в общем случае,  случайная величина, на которую можно перенести утверждения, аналогичные имеющимся в классическом интегральном исчислении. Равенство-сопоставление \eqref{Ayd2.3} в указанном смысле -- аналог замены переменных в теории
интегрирования в математическом анализе, в котором якобиан преобразования $\textit{J}_{i,j}(t)$ -- случайная величина, составленная из решений системы СДУ
(\ref{Ayd2.2}).

Определим, при каких условиях можно рассматривать вопрос об
инвариантности (сохранении) стохастического объема.

\section{Стохастическое ядро стохастического интегрального инварианта}

Рассмотрим систему обобщенных стохастических дифференциальных
уравнений Ито вида  (\ref{Ayd01.1}).
Пусть $\rho(t;{\bf x};\omega)$ -- случайная функция, измеримая
относительно потока $\sigma$-алгебр
$\Bigl\{\mathcal{F}\Bigr\}_{0}^{T}$,
$\mathcal{F}_{t_{1}}\subset\mathcal{F}_{t_{2}}$, $t_{1}<t_{2}$,
согласованного с процессами ${\bf w} (t)$ и $\nu(t;\Delta\gamma)$
(далее параметр $\omega$ будем опускать) и  относительно любой функции $f(t;{\bf x})$ из класса локально ограниченных функций, имеющей ограниченные вторые
производные по ${\bf x}$ для нее выполнены
соотношения:
\begin{equation}\label{Aydyad1}
\displaystyle\int\limits_{\mathbb{R}^{n}}\rho(t;{\bf x})f(t;{\bf
x})d\Gamma({\bf x})=\int\limits_{\mathbb{R}^{n}}\rho(0;{\bf
y})f(t;{\bf x}(t;{\bf y}))d\Gamma({\bf y})
\end{equation}
\begin{equation}\label{Aydyad2}
\displaystyle\int\limits_{\mathbb{R}^{n}}\rho(0;{\bf
x})d\Gamma({\bf x})=1,
\end{equation}
\begin{equation}\label{Aydyad3}
\begin{array}{c}
\displaystyle  \lim\limits_{|{\bf x}|\to \infty}\rho(0;{\bf x})=0,
\ \ \ d\Gamma({\bf x})=\prod\limits_{i=1}^{n}dx_{i},
\end{array}
\end{equation}
где ${\bf x}(t;{\bf y})$ -- решение уравнения (\ref{Ayd01.1}).

Если $f(t;{\bf x})=1$,  то из условия (\ref{Aydyad1}) и
(\ref{Aydyad2}) следует, что
\begin{equation}\label{usl-inv}
\displaystyle\int\limits_{\mathbb{R}^{n}}\rho(t;{\bf
x})d\Gamma({\bf x})=\int\limits_{\mathbb{R}^{n}}\rho(0;{\bf
y})d\Gamma({\bf y})=1,
\end{equation}
т. е. для  случайной функции $\rho(t;{\bf x})=\rho(t;{\bf
x};\omega)$ существует случайный функционал, сохраняющий постоянное значение:
\begin{equation}\label{ro1}
\displaystyle\int\limits_{\mathbb{R}^{n}}\rho(t;{\bf
x})d\Gamma({\bf x})=1,
\end{equation}
который можно рассматривать как стохастический объем.

Тогда \eqref{Aydyad1} с условиями \eqref{Aydyad2} и \eqref{Aydyad3} можно рассматривать как стохастический интегральный инвариант для функции $f(t;{\bf x})$.

\begin{dfn} Неотрицательную функцию $\rho(t;{\bf x})$ будем называть
стохастическим ядром или стохастической плотностью стохастического интегрального
инварианта, если выполняются равенства {\rm(\ref{Aydyad1})}, \eqref{Aydyad2} и \eqref{Aydyad3}.
\end{dfn}

\begin{remark}
Понятие интегрального инварианта для системы обыкновенных дифференциальных уравнений было известно ранее, например, оно рамматривалось В.\,И.~Зубовым в {\rm\cite[\S 8
]{Zubov_82}}. Однако существенное отличие, позволившее в {\rm\cite{D_89,D_02}} рассматривать инвариантность случайного объема на основе ядра интегрального инварианта, состоит в том, что в \eqref{Aydyad1} присутствует функциональный множитель. Таким образом, понятие ядра интегрального инварианта для системы обыкновенных дифференциальных уравнений, можно рассматривать, как частный случай введенного в {\rm\cite{D_89,D_02}}, если взять $f(t;{\bf x})=1$ и, исключив случайность, заданную винеровскими и пуассоновскими процессами, рассматривать интегрирование по детерминированному объему.
\end{remark}

Определим соотношения, при которых функция $\rho(t;{\bf x})$ для
произвольной дважды дифференцируемой функции $f(t;{\bf x})$ будет
ядром интегрального инварианта.

Для случайной функции $f(t;{\bf x}(t))$, где ${\bf x}(t)$ --
решение уравнения (\ref{Ayd01.1}), запишем обобщенную формулу Ито
{\rm{\cite[с.\, 271-272]{GS_68}}}:
\begin{equation}\label{Ayd2}
\begin{array}{c}
 \displaystyle d_{t}f(t;{\bf x}(t))=\Bigl[ \frac{\partial f(t;{\bf x}(t)) }{\partial t}+
 \sum\limits_{i=1}^{n}a_{i}(t;{\bf x}(t))\frac{\partial f(t;{\bf x}(t)) }{\partial x_{i}}+\Bigr.\\
\Bigl.+ \displaystyle  \frac{1}{2} \sum\limits_{i=1}^{n}
 \sum\limits_{j=1}^{n}
 \sum\limits_{k=1}^{m}b_{i\,k}(t;{\bf x}(t))b_{j\,k}(t;{\bf x}(t))
 \frac{\partial^{\,2} f(t;{\bf x}(t)) }{\partial x_{i}x_{j}}\Bigr]dt +\\
 +\displaystyle \sum\limits_{i=1}^{n}
 \sum\limits_{k=1}^{m}b_{i\,k}(t;{\bf x}(t))\frac{\partial f(t;{\bf x}(t)) }{\partial x_{i}} dw_{k}(t)+\\
 +\displaystyle \int\limits_{\mathbb{R}(\gamma)}
 \Bigl[ f\left(t;{\bf x}(t)+g(t;{\bf x}(t);\gamma)\right)- f(t;{\bf x}(t))
 \Bigr]\nu(dt;d\gamma).
 \end{array}
\end{equation}

\begin{remark}
В дальнейшем, для упрощения записей, будем иметь в
виду, что по индексам, встречающимся дважды, проводится
суммирование (без использования знака суммы).
\end{remark}

Продифференцируем по $t$ обе части равенства (\ref{Aydyad1}),
учитывая, что в левой части $f(t;{\bf x})$  -- детерминированная
функция, а в правой $f(t;{\bf x}(t;{\bf y}))$ -- случайный
процесс. Имеем:
\begin{equation*}
\begin{array}{c}
  \displaystyle\int\limits_{\mathbb{R}^{n}}\Bigl(f(t;{\bf
x})d_{t}\rho(t;{\bf x})+\rho(t;{\bf x})\frac{\partial f(t;{\bf
x})}{\partial t}dt\Bigr)d\Gamma({\bf
x})= \int\limits_{\mathbb{R}^{n}}\rho(0;{\bf y})d_{t}f(t;{\bf x}(t;{\bf
y}))d\Gamma({\bf y}).
\end{array}
\end{equation*}
Тогда, в силу (\ref{Aydyad1}) и (\ref{Ayd2}), получаем:

\begin{equation}\label{Aydyad5}
\begin{array}{c}
  \displaystyle\int\limits_{\mathbb{R}^{n}}\Bigl(f(t;{\bf
x})d_{t}\rho(t;{\bf x})+\rho(t;{\bf x})\frac{\partial f(t;{\bf
x})}{\partial t}dt\Bigr)d\Gamma({\bf
x})=\\
=\displaystyle\int\limits_{\mathbb{R}^{n}}\rho(0;{\bf y})d_{t}
f(t;{\bf x}(t;{\bf y}))d\Gamma({\bf y})
=\displaystyle\int\limits_{\mathbb{R}^{n}}\rho(t;{\bf
x})d_{t}f(t;{\bf
x})d\Gamma({\bf x})= \\
=\displaystyle\int\limits_{\mathbb{R}^{n}}d\Gamma({\bf
x})\rho(t;{\bf x})\cdot \biggl[ \Bigl[\frac{\partial f(t;{\bf
x}(t)) }{\partial t}+
 a_{i}(t;{\bf x})\frac{\partial f(t;{\bf x}) }{\partial x_{i}}+\Bigr.\biggr.\\
\Bigl.+ \displaystyle  \frac{1}{2} b_{i\,k}(t;{\bf
x})b_{j\,k}(t;{\bf x})
 \frac{\partial^{\,2} f(t;{\bf x}) }{\partial x_{i} \partial x_{j}}\Bigr]dt
 +\displaystyle
 b_{i\,k}(t;{\bf x})\frac{\partial f(t;{\bf x}) }{\partial x_{i}} dw_{k}(t)+\\
\biggl. +\displaystyle \int\limits_{\mathbb{R}(\gamma)}
 \Bigl[ f\left(t;{\bf x}+g(t;{\bf x};\gamma)\right)- f(t;{\bf x})
 \Bigr]\nu(dt;d\gamma)\biggr].
\end{array}
\end{equation}
Рассмотрим интеграл
\begin{equation}\label{Aydy1}
I=\displaystyle\int\limits_{\mathbb{R}^{n}}d\Gamma({\bf
x})\rho(t;{\bf x})
 f\left(t;{\bf x}+g(t;{\bf x};\gamma)\right).
\end{equation}
Сделаем замену переменных:
\begin{equation}\label{Aydx1y}
{\bf x}+g(t;{\bf x};\gamma)={\bf y},
\end{equation}
тогда, учитывая переход от ${\bf x}$ к ${\bf y}$, получаем:
$
{\bf x}={\bf y}-g(t;{\bf x};\gamma)={\bf y}-g(t;{\bf
x}^{-1}(t;{\bf y};\gamma);\gamma).
$
Следовательно, интеграл $I$ примет вид:
\begin{equation}\label{Aydy2}
\begin{array}{c}
  I=\displaystyle \int\limits_{\mathbb{R}^{n}}d\Gamma({\bf
y})\rho\left(t;{\bf y}-g(t;{\bf x}^{-1}(t;{\bf
y};\gamma);\gamma)\right)
  \cdot f(t;{\bf y})\cdot D\left( {\bf x}^{-1}(t;{\bf y};\gamma)  \right),
\end{array}
\end{equation}
где $D\left( {\bf x}^{-1}(t;{\bf y};\gamma)\right) $ -- якобиан
перехода от ${\bf x}$ к ${\bf y}$.

С учетом (\ref{Aydy2}) и интегрирования по всему пространству
$\mathbb{R}^{n}$, что дает возможность формальной замены
обозначения переменной интегрирования, получаем:
\begin{equation}\label{Aydy3}
\begin{array}{c}
  \displaystyle\int\limits_{\mathbb{R}^{n}}d\Gamma({\bf
x})\rho(t;{\bf x})
 \displaystyle \int\limits_{\mathbb{R}(\gamma)}
 \Bigl[ f\left(t;{\bf x}+g(t;{\bf x};\gamma)\right)- f(t;{\bf x})
 \Bigr]\nu(dt;d\gamma)= \\
  =\displaystyle\int\limits_{\mathbb{R}^{n}}d\Gamma({\bf
x})\rho(t;{\bf x})
 \displaystyle \int\limits_{\mathbb{R}(\gamma)}
f\left(t;{\bf x}+g(t;{\bf
x};\gamma)\right)\nu(dt;d\gamma)-\\
-\displaystyle\int\limits_{\mathbb{R}^{n}}d\Gamma({\bf
x})\rho(t;{\bf x})
 \displaystyle \int\limits_{\mathbb{R}(\gamma)} f(t;{\bf
x})
\nu(dt;d\gamma)= \\
=\displaystyle \int\limits_{\mathbb{R}^{n}}d\Gamma({\bf x})
  \int\limits_{\mathbb{R}(\gamma)}  \Bigl(\Bigr.\rho\left(t;{\bf x}-g(t;{\bf x}^{-1}(t;{\bf
x};\gamma);\gamma)\right) \cdot
 f(t;{\bf x})\cdot D\left( {\bf x}^{-1}(t;{\bf x};\gamma)  \right)\Bigl.\Bigr) \nu(dt;d\gamma)-\\
 -\displaystyle\int\limits_{\mathbb{R}^{n}}d\Gamma({\bf
x})\rho(t;{\bf x})
 \displaystyle \int\limits_{\mathbb{R}(\gamma)} f(t;{\bf x}) \nu(dt;d\gamma)=\\
=\displaystyle \int\limits_{\mathbb{R}^{n}}d\Gamma({\bf
x})f(t;{\bf x})\int\limits_{\mathbb{R}(\gamma)}
\Bigl[\rho\left(t;{\bf x}-g(t;{\bf x}^{-1}(t;{\bf
x};\gamma);\gamma)\right)
 \cdot \Bigr.\\
 \cdot\Bigl. D\left( {\bf x}^{-1}(t;{\bf x};\gamma)  \right)-\rho(t;{\bf
 x})\Bigr]\nu(dt;d\gamma)
\end{array}
\end{equation}

Учитывая (\ref{Aydyad2}), вычислим следующие интегралы, используя
интегрирование по частям:
\begin{equation}\label{Aydi01}
\displaystyle\int\limits_{\mathbb{R}^{n}}d\Gamma({\bf
x})\rho(t;{\bf x})a_{i}(t;{\bf x})\frac{\partial f(t;{\bf x})
}{\partial x_{i}},
\end{equation}
\begin{equation}\label{Aydi02}
\displaystyle\int\limits_{\mathbb{R}^{n}}d\Gamma({\bf
x})\rho(t;{\bf x})b_{i\,k}(t;{\bf x})\frac{\partial f(t;{\bf x})
}{\partial x_{i}},
\end{equation}
\begin{equation}\label{Aydi03}
\displaystyle\int\limits_{\mathbb{R}^{n}}d\Gamma({\bf
x})\rho(t;{\bf x})b_{i\,k}(t;{\bf x})b_{j\,k}(t;{\bf x})
 \frac{\partial^{\,2} f(t;{\bf x}) }{\partial x_{i} \partial x_{j}}.
\end{equation}
Для (\ref{Aydi01}) имеем:
\begin{equation*}
\begin{array}{c}
  \displaystyle\int\limits_{\mathbb{R}^{n}}d\Gamma({\bf
x})\rho(t;{\bf x})a_{i}(t;{\bf x})\frac{\partial f(t;{\bf x})
}{\partial x_{i}}
=\displaystyle\int\limits_{-\infty}^{+\infty}\prod\limits_{j=1}^{n}dx_{j}\rho(t;{\bf
x})a_{i}(t;{\bf x})\frac{\partial f(t;{\bf x}) }{\partial x_{i}}= \\
  =\displaystyle\int\limits_{-\infty}^{+\infty}\prod\limits_{j=1}^{n-1}dx_{j}
  \int\limits_{-\infty}^{+\infty}\rho(t;{\bf
x})a_{i}(t;{\bf x})\frac{\partial f(t;{\bf x}) }{\partial
x_{i}}dx_{i}.
\end{array}
\end{equation*}
С учетом (\ref{Aydyad3}) вычислим внутренний интеграл с помощью интегрирования по частям:
\begin{equation*}
\begin{array}{c}
\displaystyle\int\limits_{-\infty}^{+\infty}\rho(t;{\bf
x})a_{i}(t;{\bf x})\frac{\partial f(t;{\bf x}) }{\partial
x_{i}}dx_{i}=\\
=\rho(t;{\bf x})a_{i}(t;{\bf x})f(t;{\bf
x})\biggl|_{-\infty}^{+\infty}-\displaystyle\int\limits_{-\infty}^{+\infty}f(t;{\bf
x})\displaystyle\frac{\partial \left( \rho(t;{\bf x})a_{i}(t;{\bf
x})\right) }{\partial
x_{i}}dx_{i} =\\
=-\displaystyle\int\limits_{-\infty}^{+\infty}f(t;{\bf
x})\displaystyle\frac{\partial \left( \rho(t;{\bf x})a_{i}(t;{\bf
x})\right) }{\partial x_{i}}dx_{i} .
\end{array}
\end{equation*}
Следовательно,
\begin{equation}\label{Aydi10}
\begin{array}{c}
\displaystyle\int\limits_{\mathbb{R}^{n}}d\Gamma({\bf
x})\rho(t;{\bf x})a_{i}(t;{\bf x})\frac{\partial f(t;{\bf x})
}{\partial x_{i}}
=\displaystyle\int\limits_{-\infty}^{+\infty}\prod\limits_{j=1}^{n}dx_{j}\rho(t;{\bf
x})a_{i}(t;{\bf x})\frac{\partial f(t;{\bf x}) }{\partial x_{i}}= \\
  =\displaystyle\int\limits_{-\infty}^{+\infty}\prod\limits_{j=1}^{n-1}dx_{j}
  \int\limits_{-\infty}^{+\infty}\rho(t;{\bf
x})a_{i}(t;{\bf x})\frac{\partial f(t;{\bf x}) }{\partial
x_{i}}dx_{i}=\\
=-\displaystyle\int\limits_{-\infty}^{+\infty}\prod\limits_{j=1}^{n}dx_{j}
\displaystyle\int\limits_{-\infty}^{+\infty}f(t;{\bf
x})\displaystyle\frac{\partial \left( \rho(t;{\bf x})a_{i}(t;{\bf
x})\right) }{\partial
x_{i}}dx_{i} =
-\displaystyle\int\limits_{\mathbb{R}^{n}}d\Gamma({\bf
x})f(t;{\bf x})\displaystyle\frac{\partial \left( \rho(t;{\bf
x})a_{i}(t;{\bf x})\right) }{\partial x_{i}}
\end{array}
\end{equation}
Аналогичным образом вычислим второй интеграл (\ref{Aydi02})  и,
применяя дважды интегрирование по частям, вычислим третий интеграл
(\ref{Aydi03}):
\begin{equation}\label{Aydi20}
\begin{array}{c}
  \displaystyle\int\limits_{\mathbb{R}^{n}}d\Gamma({\bf
x})\rho(t;{\bf x})b_{i\,k}(t;{\bf x})\frac{\partial f(t;{\bf x})
}{\partial x_{i}}
  =-\displaystyle\int\limits_{\mathbb{R}^{n}}d\Gamma({\bf
x})f(t;{\bf x})\frac{\partial\rho(t;{\bf x})b_{i\,k}(t;{\bf x})
}{\partial x_{i}},
\end{array}
\end{equation}
\begin{equation}\label{Aydi30}
\begin{array}{c}
  \displaystyle\int\limits_{\mathbb{R}^{n}}d\Gamma({\bf
x})\rho(t;{\bf x})b_{i\,k}(t;{\bf x})b_{j\,k}(t;{\bf x})
 \frac{\partial^{\,2} f(t;{\bf x}) }{\partial x_{i} \partial x_{j}}=\\
 =\displaystyle\int\limits_{\mathbb{R}^{n}}d\Gamma({\bf x})f(t;{\bf
x})
 \frac{\partial^{\,2} \left(\rho(t;{\bf x})b_{i\,k}(t;{\bf x})b_{j\,k}(t;{\bf x})\right) }{\partial x_{i} \partial
 x_{j}}.
\end{array}
\end{equation}
В (\ref{Aydyad5}) перенесем все в правую часть и с учетом
(\ref{Aydy3}),
 (\ref{Aydi10}), (\ref{Aydi20}) и (\ref{Aydi30}) получим:
\begin{equation}\label{Aydi}
\begin{array}{c}
 0= \displaystyle\int\limits_{\mathbb{R}^{n}}d\Gamma({\bf x})f(t;{\bf
x})\cdot \biggl[-d_{t}\left(\rho(t;{\bf x})\right)-\frac{\partial
f(t;{\bf x}(t)) }{\partial t}+\biggr.\\
+ \displaystyle\frac{\partial f(t;{\bf x}(t)) }{\partial
t}-\frac{\partial\rho(t;{\bf x})b_{i\,k}(t;{\bf x}) }{\partial
x_{i}}dw_{k}(t)+\\
+\Bigl(\Bigr.-\displaystyle\frac{\partial \left( \rho(t;{\bf
x})a_{i}(t;{\bf x}(t))\right) }{\partial x_{i}}
+\displaystyle\frac{1}{2}\frac{\partial^{\,2} \left(\rho(t;{\bf
x})b_{i\,k}(t;{\bf x})b_{j\,k}(t;{\bf x})\right) }{\partial x_{i} \partial x_{j}}\Bigl.\Bigr)dt + \\
 \biggl.+\displaystyle\int\limits_{\mathbb{R}(\gamma)} \Bigl[\rho\left(t;{\bf
x}-g(t;{\bf x}^{-1}(t;{\bf x};\gamma);\gamma)\right)
  \cdot D\left( {\bf x}^{-1}(t;{\bf x};\gamma)  \right)-\rho(t;{\bf
 x})\Bigr]\nu(dt;d\gamma)\biggr].
\end{array}
\end{equation}
Чтобы равенство (\ref{Aydi}) имело место для любой локально
ограниченной функции $f(t;{\bf x})$,  имеющей ограниченные
производные второго порядка, интегральный инвариант $\rho(t;{\bf
x})$ должен являться решением стохастического уравнения
\begin{equation}\label{Aydii}
\begin{array}{c}
 d_{t}\rho(t;{\bf x})=
-\displaystyle\frac{\partial\rho(t;{\bf x})b_{i\,k}(t;{\bf x})
}{\partial x_{i}}dw_{k}(t) +\Bigl(-\displaystyle\frac{\partial
\left( \rho(t;{\bf x})a_{i}(t;{\bf x})\right) }{\partial
x_{i}}+\\
+\displaystyle\frac{1}{2}\frac{\partial^{\,2} \left(\rho(t;{\bf
x})b_{i\,k}(t;{\bf x})b_{j\,k}(t;{\bf x})\right) }{\partial x_{i} \partial x_{j}}\Bigr)dt + \\
 +\displaystyle\int\limits_{\mathbb{R}(\gamma)} \Bigl[\rho\left(t;{\bf
x}-g(t;{\bf x}^{-1}(t;{\bf x};\gamma);\gamma)\right)
 \cdot D\left( {\bf x}^{-1}(t;{\bf x};\gamma)  \right)-\rho(t;{\bf
 x})\Bigr]\nu(dt;d\gamma).
\end{array}
\end{equation}
При этом должны выполняться условия
\begin{equation}\label{Aydii1}
\begin{array}{c}
\rho(t;{\bf x})\Bigl|_{t=0}=\rho(0;{\bf x})\in C_{0}^{2},  \\
\lim\limits_{|{\bf x}|\to \infty}\rho(0;{\bf x})=0, \ \ \ \
\displaystyle\lim\limits_{|{\bf x}|\to
\infty}\frac{\partial\rho(0;{\bf x})}{\partial x_{i}}=0.
 \end{array}
\end{equation}
Таким образом, получены условия инвариантности стохастического объема и доказана следующая
\begin{theorem}\label{thro1}
Пусть ${\bf x}(t)$,
${\bf x}\in\mathbb{R}^{n}$, решение системы обобщенных стохастических
дифференциальных уравнений Ито
\begin{equation*}
\begin{array}{c}
  dx_{i}(t)= \displaystyle a_{i}\Bigl(t;{\bf x}(t) \Bigr)\, dt
  +\sum\limits_{k=1}^{m}
b_{i,k}\Bigl(t;{\bf x}(t)\Bigr)\, dw_{k}(t) +
\int\limits_{\mathbb{R}(\gamma)}g_{i}\Bigl(t;{\bf x}(t);\gamma\Bigr)\nu(dt;d\gamma), \\
  {\bf x}(t)={\bf x}\Bigl(t;{\bf x}(0) \Bigr) \Bigr|_{t=0}={\bf
  x}(0), \ \ \ \ \ \ i=\overline{1,n}, \ \ \ \ \ \ t\geq 0,
\end{array}
\end{equation*}
где ${\bf w} (t)$ -- $m$-мерный винеровский процесс,
$\nu(t;\Delta\gamma)$ -- однородная по $t$ нецентрированная мера
Пуассона и $\rho(t;{\bf x})$ -- случайная функция, измеримая
относительно потока $\sigma$-алгебр
$\Bigl\{\mathcal{F}\Bigr\}_{0}^{T}$,
$\mathcal{F}_{t_{1}}\subset\mathcal{F}_{t_{2}}$, $t_{1}<t_{2}$,
согласованного с процессами ${\bf w} (t)$ и $\nu(t;\Delta\gamma)$
и  относительно любой функции $f(t;{\bf x})$ из класса $\mathfrak{S}$ локально ограниченных функций, имеющих ограниченные вторые
производные по ${\bf x}$. Функция $\rho(t;{\bf x})$ является стохастическом ядром стохастического интегрального инварианта
\eqref{Aydyad1} для
произвольной локально ограниченной функции $f(t;{\bf x})\in \mathfrak{S}$, если она является решением системы \eqref{Aydii} обобщенных СДУ Ито, удовлетворяющим начальным условиям \eqref{Aydii1}.
\end{theorem}

\section{Обобщение формулы Ито-Вентцеля}

Равенство (\ref{Aydyad1}) с условиями (\ref{Aydyad2}), (\ref{Aydyad3}) и уравнение (\ref{Aydii}) для ядра
интегрального инварианта соответствовали случаю детерминированной
функции $f(t;{\bf x})$.  Положим выполнение аналогичного равенства  для случайной функции ${\bf
z}(t;{\bf x})$:
\begin{equation}\label{Aydz3p}
\begin{array}{c}
  \displaystyle\int\limits_{\mathbb{R}^{n}}\rho(t;{\bf x}){\bf z}(t;{\bf
x})d\Gamma({\bf x})=
  \displaystyle\int\limits_{\mathbb{R}^{n}}\rho(0;{\bf y}){\bf
z}(t;{\bf x}(t;{\bf y}))d\Gamma({\bf y}),
\end{array}
\end{equation}
где ${\bf x}(t;{\bf y})$ -- решение системы СДУ (\ref{Ayd01.1}).

Рассмотрим сложный случайный процесс ${\bf z}\left(t;{\bf x}(t;{\bf
y})\right)\in \mathbb{R}^{n_{o}}$, где ${\bf x}(t;{\bf y})$ --
решение системы СДУ (\ref{Ayd01.1}), а процесс ${\bf z}(t;{\bf
x})$ -- решение системы обобщенных СДУ Ито:
\begin{equation}\label{Aydz1}
\begin{array}{c}
  d_{t}{\bf z}(t;{\bf x})=\Pi(t;{\bf x})dt+D_{k}(t;{\bf x})dw_{k}(t)
+
  \displaystyle\int\limits_{R(\gamma)}\nu(dt;d\gamma)G(t;{\bf x};\gamma).
  \end{array}
\end{equation}
Относительно случайных функций $\Pi(t;{\bf x})$, $D_{k}(t;{\bf
x})$, $G(t;{\bf x};\gamma)$, определенных на пространстве
$\mathbb{R}^{n_{o}}$, предполагаем, что они непрерывны и
ограничены вместе со своими производными по всем переменным,
измеримые относительно потока $\sigma$-алгебр $\mathcal{F}_{t}$,
согласованного с процессами ${\bf w}(t)$ и $\nu(t;\Delta\gamma)$
из (\ref{Ayd01.1}).

Опираясь на уравнение для стохастического интегрального
инварианта, построим правило дифференцирования
для сложного случайного процесса ${\bf z}\left(t;{\bf x}(t;{\bf
y})\right)$.

Рассмотрим интеграл
$
\displaystyle\int\limits_{\mathbb{R}^{n}}\rho(t;{\bf x}){\bf
z}(t;{\bf x})d\Gamma({\bf x})
$.
Поскольку интегрирование проводится по пространству
$\mathbb{R}^{n}$, на котором задан процесс ${\bf x}(t)$,
используем равенство (\ref{Aydz3p}), записав его в виде (поменяв части равенства местами):
\begin{equation}\label{Aydz3}
\begin{array}{c}
  \displaystyle\int\limits_{\mathbb{R}^{n}}\rho(0;{\bf y}){\bf
z}(t;{\bf x}(t;{\bf y}))d\Gamma({\bf y})=
  \int\limits_{\mathbb{R}^{n}}\rho(t;{\bf x}){\bf z}(t;{\bf
x})d\Gamma({\bf x}).
\end{array}
\end{equation}
Продифференцируем обе части (\ref{Aydz3}) по $t$:
\begin{equation}\label{Aydz4}
\begin{array}{c}
  \displaystyle\int\limits_{\mathbb{R}^{n}}\rho(0;{\bf y})d_{t}{\bf
z}(t;{\bf x}(t;{\bf y}))d\Gamma({\bf
y})= \\
 \displaystyle\int\limits_{\mathbb{R}^{n}}\Bigl(\Bigr.\rho(t;{\bf x})d_{t}{\bf
z}(t;{\bf x})+{\bf z}(t;{\bf x})d_{t}\rho(t;{\bf x}) -
D_{k}(t;{\bf x}) \displaystyle\frac{\partial\rho(t;{\bf
x})b_{i\,k}(t;{\bf x}) }{\partial x_{i}}dt\Bigl.\Bigr)d\Gamma({\bf
x}) ,
\end{array}
\end{equation}
Учитывая, что интегрирование идет по пространству $\mathbb{R}^{n}$
после введения замены переменной, то при интегрировании по кривой
$R(\gamma)$ в этом пространстве нужно учитывать произведенную
замену (\ref{Aydx1y}).

Поскольку функция $\rho(t;{\bf x})$ -- ядро интегрального инварианта \eqref{Aydz3p}, применим теорему~\ref{thro1}.  Подставим (\ref{Aydii}) и (\ref{Aydz1}) в правую часть
(\ref{Aydz4}):
$$
I_{1}= \displaystyle \int\limits_{\mathbb{R}^{n}}\Bigl(\rho(t;{\bf
x})d_{t}{\bf z}(t;{\bf x})+{\bf z}(t;{\bf x})d_{t}\rho(t;{\bf x})
-D_{k}(t;{\bf x}) \displaystyle\frac{\partial\rho(t;{\bf
x})b_{i\,k}(t;{\bf x})
}{\partial x_{i}}dt\Bigr)d\Gamma({\bf x})=
$$
$$  = \displaystyle\int\limits_{\mathbb{R}^{n}}d\Gamma({\bf x})\biggl\{\rho(t;{\bf x})
  \Bigl[\Bigr.\Pi(t;{\bf x})dt+D_{k}(t;{\bf x})dw_{k}(t)
+
$$
$$
+\displaystyle\int\limits_{R(\gamma)}G(t;{\bf x}+g(t;{\bf
x}^{-1}(t;{\bf
x};\gamma);\gamma)\nu(dt;d\gamma)\Bigl.\Bigr]+
$$
$$
+{\bf z}(t;{\bf x})\biggl[-\displaystyle\frac{\partial\rho(t;{\bf
x})b_{i\,k}(t;{\bf x}) }{\partial
x_{i}}dw_{k}(t)+\Bigr.\displaystyle\Bigl(\Bigr.-\displaystyle\frac{\partial
\left( \rho(t;{\bf x})a_{i}(t;{\bf x})\right) }{\partial
x_{i}}+
\displaystyle\frac{1}{2}\frac{\partial^{\,2} \left(\rho(t;{\bf
x})b_{i\,k}(t;{\bf x})b_{j\,k}(t;{\bf x})\right) }{\partial x_{i} \partial x_{j}}\Bigl.\Bigr)dt +
$$
$$
 +\displaystyle\int\limits_{\mathbb{R}(\gamma)} \Bigl[\Bigr.\rho\left(t;{\bf
x}-g(t;{\bf x}^{-1}(t;{\bf x};\gamma);\gamma)\right)
 \cdot  D\left( {\bf x}^{-1}(t;{\bf x};\gamma)  \right)-\rho(t;{\bf
 x})\Bigl.\Bigr]\nu(dt;d\gamma)\Bigl.\Bigl.\biggr]-
 $$
 $$
 -D_{k}(t;{\bf
x}) \displaystyle\frac{\partial\rho(t;{\bf x})b_{i\,k}(t;{\bf x})
}{\partial x_{i}}dt\biggr\}.
$$
Приведем подобные:
\begin{equation}\label{Aydz5}
\begin{array}{c}
I_{1}=\displaystyle\int\limits_{\mathbb{R}^{n}}d\Gamma({\bf
x})\biggl(\rho(t;{\bf x})\Pi(t;{\bf x}) -{\bf z}(t;{\bf
x})\displaystyle\frac{\partial \left( \rho(t;{\bf x})a_{i}(t;{\bf
x})\right) }{\partial x_{i}}+\biggr.\\
-D_{k}(t;{\bf x}) \displaystyle\frac{\partial\rho(t;{\bf
x})b_{i\,k}(t;{\bf x}) }{\partial x_{i}}
+\displaystyle\frac{1}{2}{\bf z}(t;{\bf x})\frac{\partial^{\,2}
\left(\rho(t;{\bf x})b_{i\,k}(t;{\bf x})b_{j\,k}(t;{\bf x})\right)
}{\partial x_{i}
\partial x_{j}}\biggl.\biggr)dt+\\
+\displaystyle\int\limits_{\mathbb{R}^{n}}d\Gamma({\bf
x})\biggl(\rho(t;{\bf x})D_{k}(t;{\bf x}) -{\bf z}(t;{\bf
x})\displaystyle\frac{\partial \left( \rho(t;{\bf
x})b_{i\,k}(t;{\bf x})\right) }{\partial
x_{i}}\biggr)dw_{k}(t)+\\
+\displaystyle\int\limits_{\mathbb{R}^{n}}d\Gamma({\bf
x})\biggl[\rho(t;{\bf x}) \cdot\int\limits_{R(\gamma)}G(t;{\bf
x}+g(t;{\bf x}^{-1}(t;{\bf
x};\gamma);\gamma))\nu(dt;d\gamma) +\Bigr.\\
+{\bf z}(t;{\bf
x})\displaystyle\int\limits_{R(\gamma)}\Bigl[\rho\left(t;{\bf
x}-g(t;{\bf x}^{-1}(t;{\bf x};\gamma);\gamma)\right)
 \cdot D\left( {\bf x}^{-1}(t;{\bf x};\gamma)  \right)-\rho(t;{\bf
 x})\Bigr]\nu(dt;d\gamma)\Bigl.\biggr].
\end{array}
\end{equation}
В силу (\ref{Aydi10}), (\ref{Aydi20}), (\ref{Aydi30}) имеем:
\begin{equation*}
\begin{array}{c}
\displaystyle\int\limits_{\mathbb{R}^{n}}d\Gamma({\bf x}){\bf
z}(t;{\bf x})\frac{\partial \left(\rho(t;{\bf x})a_{i}(t;{\bf
x})\right) }{\partial x_{i}} =
-\displaystyle\int\limits_{\mathbb{R}^{n}}d\Gamma({\bf
x})\rho(t;{\bf x})a_{i}(t;{\bf x})\displaystyle\frac{\partial
 {\bf z}(t;{\bf x})
 }{\partial x_{i}},
\end{array}
\end{equation*}
\begin{equation*}
\begin{array}{c}
  \displaystyle\int\limits_{\mathbb{R}^{n}}d\Gamma({\bf
x})D_{k}(t;{\bf x})\frac{\partial\left(\rho(t;{\bf
x})b_{i\,k}(t;{\bf x})\right) }{\partial x_{i}}=
  -\displaystyle\int\limits_{\mathbb{R}^{n}}d\Gamma({\bf
x})\rho(t;{\bf x})b_{i\,k}(t;{\bf x})\frac{\partial D_{k}(t;{\bf
x}) }{\partial x_{i}},
\end{array}
\end{equation*}
\begin{equation*}
\begin{array}{c}
\displaystyle\int\limits_{\mathbb{R}^{n}}d\Gamma({\bf x}){\bf
z}(t;{\bf x})
 \frac{\partial^{\,2} \left(\rho(t;{\bf x})b_{i\,k}(t;{\bf x})b_{j\,k}(t;{\bf x})\right) }{\partial x_{i} \partial
 x_{j}}=
  \\
 = \displaystyle\int\limits_{\mathbb{R}^{n}}d\Gamma({\bf
x})\rho(t;{\bf x})b_{i\,k}(t;{\bf x})b_{j\,k}(t;{\bf x})
 \frac{\partial^{\,2} {\bf z}(t;{\bf x}) }{\partial x_{i} \partial
 x_{j}}.
\end{array}
\end{equation*}
Вычислим последний интеграл в сумме (\ref{Aydz5}):
\begin{equation}\label{Ayde6}
\begin{array}{c}
  I_{2}=\displaystyle\int\limits_{\mathbb{R}^{n}}d\Gamma({\bf
x}){\bf z}(t;{\bf x})
  \cdot\rho\Bigl(t;{\bf x}-g(t;{\bf x}^{-1}(t;{\bf
x};\gamma);\gamma)\Bigr)D({\bf x}^{-1}(t;{\bf x};\gamma)).
\end{array}
\end{equation}
Введем замену переменных:
$$
\begin{array}{c}
  {\bf x}-g(t;{\bf x}^{-1}(t;{\bf x};\gamma);\gamma)={\bf y}, \\
  {\bf x}={\bf y}+g(t;{\bf x}^{-1}(t;{\bf x};\gamma);\gamma)={\bf y}+g(t;{\bf
  y};\gamma).
\end{array}
$$
Обозначим якобиан перехода от вектора ${\bf x}$ к вектору ${\bf
y}$ через $D_{o}({\bf x}^{-1}(t;{\bf y};\gamma))$. Тогда, в силу
(\ref{Aydy1}) и дальнейшей формальной замены обозначения
переменной интегрирования, получаем:
\begin{equation}\label{Ayde7}
\begin{array}{c}
  I_{2}=\displaystyle\int\limits_{\mathbb{R}^{n}}d\Gamma({\bf y}){\bf
z}\Bigl(t;{\bf y}+g(t;{\bf
  y};\gamma)\Bigr)
  \cdot\rho(t;{\bf y})D_{o}({\bf x}^{-1}(t;{\bf y};\gamma))D({\bf
x}^{-1}(t;{\bf x};\gamma))= \\
 =\displaystyle\int\limits_{\mathbb{R}^{n}}d\Gamma({\bf y}){\bf
z}\Bigl(t;{\bf y}+g(t;{\bf
  y};\gamma)\Bigr)\rho(t;{\bf y})
  =\displaystyle\int\limits_{\mathbb{R}^{n}}d\Gamma({\bf x}){\bf
z}\Bigl(t;{\bf x}+g(t;{\bf
  x};\gamma)\Bigr)\rho(t;{\bf x}).
\end{array}
\end{equation}
В результате, правая часть (\ref{Aydz4}) примет вид:
\begin{equation}\label{Ayde8}
\begin{array}{c}
I_{1}=\displaystyle\int\limits_{\mathbb{R}^{n}}d\Gamma({\bf
x})\rho(t;{\bf x})\biggl\{\Bigl( D_{k}(t;{\bf x}) +b_{i\,k}(t;{\bf
x})\frac{\partial {\bf z}(t;{\bf x}) }{\partial x_{i}}
\Bigr)dw_{k}(t)+\biggr.\\+
 \Bigl(\Pi(t;{\bf x})+a_{i}(t;{\bf
x})\displaystyle\frac{\partial
 {\bf z}(t;{\bf x})
 }{\partial x_{i}}+ b_{i\,k}(t;{\bf
x})\frac{\partial
 D_{k}(t;{\bf x})
 }{\partial x_{i}}
 +\displaystyle\frac{1}{2}b_{i\,k}(t;{\bf x})b_{j\,k}(t;{\bf x})
 \frac{\partial^{\,2} {\bf z}(t;{\bf x}) }{\partial x_{i} \partial
 x_{j}}\Bigr)dt+\\
 +\displaystyle\int\limits_{R(\gamma)}G(t;{\bf x}+g(t;{\bf
x}^{-1}(t;{\bf x};\gamma);\gamma)\nu(dt;d\gamma)+\\
+
 \displaystyle\int\limits_{R(\gamma)}\Bigl[{\bf z}\Bigl(t;{\bf x}+g(t;{\bf
  x};\gamma)\Bigr)
  -{\bf z}(t;{\bf
 x})\Bigr]\nu(dt;d\gamma)\biggl.\biggr\}.
 \end{array}
\end{equation}
В (\ref{Aydz4}) перенесем все в правую часть и с учетом
(\ref{Aydyad1}) получаем:
\begin{equation*}
\begin{array}{c}
0=\displaystyle\int\limits_{\mathbb{R}^{n}}d\Gamma({\bf
y})\rho(0;{\bf y})\biggl\{-d_{t}{\bf z}(t;{\bf x}(t;{\bf
y}))+\biggr.\Bigl( D_{k}(t;{\bf x}(t;{\bf y}))+b_{i\,k}(t;{\bf x}(t;{\bf
y}))\displaystyle\frac{\partial {\bf z}(t;{\bf
x}) }{\partial x_{i}} \Bigr)dw_{k}(t)+\\
+ \Bigl(\Pi(t;{\bf x}(t;{\bf y}))+a_{i}(t;{\bf
x})\displaystyle\frac{\partial
 {\bf z}(t;{\bf x})
 }{\partial x_{i}}
 +b_{i\,k}(t;{\bf
x})\frac{\partial
 D_{k}(t;{\bf x})
 }{\partial x_{i}}+\Bigr.\\
 +\Bigl.\displaystyle\frac{1}{2}b_{i\,k}(t;{\bf x}(t;{\bf y}))b_{j\,k}(t;{\bf x}(t;{\bf y}))
 \frac{\partial^{\,2} {\bf z}(t;{\bf x}(t;{\bf y})) }{\partial x_{i} \partial
 x_{j}}\Bigr)dt+\\
 +\displaystyle\int\limits_{R(\gamma)}G(t;{\bf x}(t;{\bf y})+g(t;{\bf
  x}(t;{\bf y});\gamma)\nu(dt;d\gamma)+\\
 +
 \displaystyle\int\limits_{R(\gamma)}\Bigl[{\bf z}\Bigl(t;{\bf x}(t;{\bf y})+g(t;{\bf
  x}(t;{\bf y});\gamma)\Bigr)-
{\bf z}(t;{\bf
 x}(t;{\bf y}))\Bigr]\nu(dt;d\gamma)\biggl.\biggr\}.
 \end{array}
\end{equation*}
Следовательно, построенное правило дифференцирования сложного процесса будет иметь
вид, который определяет следующая
\begin{theorem}{\rm\cite{11_KchOboz}}\label{thro2}
Пусть случайный процесс ${\bf z}\left(t;{\bf x}(t;{\bf
y})\right)\in \mathbb{R}^{n_{o}}$, где ${\bf x}(t;{\bf y})\in\mathbb{R}^{n} $ --
решение системы обобщенных СДУ \eqref{Ayd01.1}, а процесс ${\bf z}(t;{\bf
x})$ -- решение системы обобщенных СДУ \eqref{Aydz1}.
Относительно коэффициентов-случайных функций, входящих в системы \eqref{Ayd01.1} и \eqref{Aydz1}, определенных на пространствах
$\mathbb{R}^{n_{o}}$ и $\mathbb{R}^{n}$ соответственно, предполагаем, что они
непрерывны и
ограничены вместе со своими производными по всем переменным,
измеримые относительно потока $\sigma$-алгебр $\mathcal{F}_{t}$,
согласованного с процессами ${\bf w}(t)$ и $\nu(t;\Delta\gamma)$
из \eqref{Ayd01.1}.
Тогда сложный случайный процесс ${\bf z}\left(t;{\bf x}(t;{\bf
y})\right)$ является решением системы обобщенных СДУ
\begin{equation}\label{Ayde9}
\begin{array}{c}
d_{t}{\bf z}(t;{\bf x}(t;{\bf y}))=\displaystyle\Bigl(
D_{k}(t;{\bf x}(t;{\bf y})) +b_{i\,k}(t;{\bf x}(t;{\bf
y}))\frac{\partial {\bf z}(t;{\bf
x}(t;{\bf y})) }{\partial x_{i}} \Bigr)dw_{k}(t)+\\
 +\Bigl(\Pi(t;{\bf x}(t;{\bf y}))+a_{i}(t;{\bf
x}(t;{\bf y}))\displaystyle\frac{\partial
 {\bf z}(t;{\bf x}(t;{\bf y}))
 }{\partial x_{i}}
 +b_{i\,k}(t;{\bf
x}(t;{\bf y}))\displaystyle\frac{\partial
 D_{k}(t;{\bf x}(t;{\bf y}))
 }{\partial x_{i}}+\Bigr.\\
 +\Bigl.\displaystyle\frac{1}{2}b_{i\,k}(t;{\bf x}(t;{\bf y}))b_{j\,k}(t;{\bf x}(t;{\bf y}))
 \frac{\partial^{\,2} {\bf z}(t;{\bf x}(t;{\bf y})) }{\partial x_{i} \partial
 x_{j}}\Bigr)dt+\\
 +\displaystyle\int\limits_{R(\gamma)}G(t;{\bf x}(t;{\bf y})+g(t;{\bf
  x}(t;{\bf y});\gamma)\nu(dt;d\gamma)+\\
 +
 \displaystyle\int\limits_{R(\gamma)}\Bigl[{\bf z}\Bigl(t;{\bf x}(t;{\bf y})+g(t;{\bf
  x}(t;{\bf y});\gamma)\Bigr)
  -{\bf z}(t;{\bf
 x}(t;{\bf y}))\Bigr]\nu(dt;d\gamma).
 \end{array}
\end{equation}
\end{theorem}
По аналогии с известной терминологией, формулу (\ref{Ayde9})
назовем обобщенной формулой Ито-Вентцеля для обобщенного уравнения
Ито с нецентрированной мерой.

Кроме того, можно сформулировать еще одно
\begin{prp}\label{thro3}
Если ${\bf z}(t;{\bf x}(t;{\bf y}))$ -- решение
системы \eqref{Ayde9}, удовлетворяющее начальному условию
\begin{equation*}
{\bf z}(t;{\bf x}(t;{\bf y}))\Bigr|_{t=0}={\bf z}(0;{\bf y}), \ \
\ {\bf z}(0;{\bf y})\in C_{o}^{2},
\end{equation*}
тогда случайная функция
$\rho(t;{\bf x})$, измеримая
относительно потока $\sigma$-алгебр
$\Bigl\{\mathcal{F}\Bigr\}_{0}^{T}$,
$\mathcal{F}_{t_{1}}\subset\mathcal{F}_{t_{2}}$, $t_{1}<t_{2}$,
согласованного с процессами ${\bf w} (t)$ и $\nu(t;\Delta\gamma)$
является стохастическом ядром стохастического интегрального инварианта
\eqref{Aydz3p}.
\end{prp}

\section{Стохастический первый интеграл}\label{PintSt}
В \cite{D_78} было введено понятие первого интеграла для системы
стохастических дифференциальных уравнений Ито (без пуассоновской
составляющей), в \cite[c.\,24]{D_02} -- понятие стохастического
первого интеграла для системы обобщенных стохастических
дифференциальных уравнений Ито с центрированной пуассоновской
мерой. Введем аналогичное понятие для случая наличия
нецентрированной меры Пуассона.

\begin{dfn}\label{Ayddf1}
Случайную функцию $u(t;{\bf x};\omega)$, определенную на
том же вероятностном пространстве, что и решение системы
\eqref{Ayd01.1}, будем называть стохастическим
первым интегралом системы \eqref{Ayd01.1} обобщенных СДУ Ито с нецентрированной пуассоновской мерой, если с
вероятностью, равной 1, выполняется условие
$$
u\Bigl(t;{\bf x}(t; {\bf x}(0));\omega\Bigr)=u\Bigl(0;{\bf
x}(0);\omega\Bigr)
$$
для любого решения ${\bf x}(t;{\bf x}(0);\omega)$ системы
\eqref{Ayd01.1}.
\end{dfn}

Определим условия, при выполнении которых  функция $u(t;{\bf x};\omega)$ будет стохастическим первым интегралом системы \eqref{Ayd01.1}.
\begin{lemma}\label{l1}
Если функция $\rho(t;{\bf x})$ -- стохастическое ядро интегрального инварианта $n$-го порядка стохастического процесса ${\bf x}(t)$, выходящего из случайной точки ${\bf x}(0)={\bf y}$, то для любого $t$ она удовлетворяет равенству  $\rho\Bigl(t;{\bf x}(t;{\bf y})\Bigr)\textit{J}(t;{\bf
y})=\rho(0;{\bf y})$, где $\textit{J}(t;{\bf
y})=\textit{J}(t;{\bf x}(0))$ -- якобиан перехода от ${\bf x}(t)$ к ${\bf x}(0)={\bf y}$.
\end{lemma}
\begin{proof}
Обратимся к равенству \eqref{usl-inv}. Произведем замену переменных в интеграле и получим утверждение леммы, поскольку интегрирование происходит по одному и тому же случайному объему:
\begin{equation*}
1=\int\limits_{\mathbb{R}^{n}}\rho(0;{\bf
y})d\Gamma({\bf y})=\displaystyle\int\limits_{\mathbb{R}^{n}}\rho(t;{\bf
x})d\Gamma({\bf x})=\displaystyle\int\limits_{\mathbb{R}^{n}}\rho\Bigl(t;{\bf x}(t;{\bf y})\Bigr)\textit{J}(t;{\bf
y})d\Gamma({\bf y}).
\end{equation*}
\end{proof}

В силу того, что случайный процесс ${\bf x}(t)$ определен в расширенном фазовом пространстве общей размерности, равной $n+1$, то для определения единственности траектории, система дифференциальных уравнений для стохастических ядер должна состоять не менее, чем из $n+1$ уравнений:
\begin{equation}\label{sys-yad}
\left\{
\begin{array}{l}
 d_{t}\rho_{l}(t;{\bf x})=
-\displaystyle\frac{\partial\rho_{l}(t;{\bf x})b_{i\,k}(t;{\bf x})
}{\partial x_{i}}dw_{k}(t) +\Bigl(-\displaystyle\frac{\partial
\left( \rho_{l}(t;{\bf x})a_{i}(t;{\bf x})\right) }{\partial
x_{i}}+\\
+\displaystyle\frac{1}{2}\frac{\partial^{\,2} \left(\rho_{l}(t;{\bf
x})b_{i\,k}(t;{\bf x})b_{j\,k}(t;{\bf x})\right) }{\partial x_{i} \partial x_{j}}\Bigr)dt + \\
 +\displaystyle\int\limits_{\mathbb{R}(\gamma)} \Bigl[\rho_{l}\left(t;{\bf
x}-g(t;{\bf x}^{-1}(t;{\bf x};\gamma);\gamma)\right)
 \cdot D\left( {\bf x}^{-1}(t;{\bf x};\gamma)  \right)-\rho_{l}(t;{\bf
 x})\Bigr]\nu(dt;d\gamma), \\
\rho_{l}(t;{\bf x}(t))\Bigl.\Bigr|_{t=0}=\rho_{l}(0;{\bf x}(0))=\rho_{l}(0;{\bf y}), \ \ \ \ \  l=\overline{1, n+1}.
\end{array}
\right.
\end{equation}

Как известно, совокупность ядер называется полной, если любая другая функция,
являющаяся ядром интегрального инварианта $n$-го порядка, может
быть представлена как функция от элементов этой совокупности.
\begin{theorem}\label{th-yad}
Система стохастических уравнений \eqref{Ayd01.1}
обладает полной со\-во\-куп\-ностью ядер, состоящей  из $(n+1)$-й  функций,
каждая из которых является решением уравнения \eqref{Aydii}.
\end{theorem}

\begin{proof} Пусть $\rho_{l}(t;{\bf x})$,  $l=\overline{1,m}$, $m\geq n+1$ --  ядра интегрального инварианта \eqref{Aydyad1}.  Из леммы \ref{l1} следует,  что
 для любых $l\neq n+1$ отношение $\dfrac{\rho_{l}(t;{\bf x}(t;{\bf y}))}
{\rho_{n+1}(t;{\bf x}(t;{\bf y}))}$ есть постоянная, зависящая только от ${\bf x}(0)={\bf y}$ для любого решения ${\bf x}(t)$ системы ДУ \eqref{Ayd01.1}:
$
\dfrac{\rho_{l}(t;{\bf x}(t;{\bf y}))}
{\rho_{n+1}(t;{\bf x}(t;{\bf y}))}=\dfrac{\rho_{l}(0;{\bf y})}
{\rho_{n+1}(0;{\bf y})}.
$
Построим функции $\theta_{l}(t;{\bf x})= \rho_{l}(t;{\bf x})\,\rho_{n+1}^{-1}(t;{\bf
x})$, $l=\overline{1,n}$,
при условии линейной независимости функций $\rho_{l}(0;{\bf y})$, $\rho_{n+1}(0;{\bf y})$.
В силу специфики подхода, описанного во введении, и полученной независимости функций $\theta_{l}(t;{\bf x})$, при $t=0$ можно было построить взаимно-однозначное соответствие:
\begin{equation}\label{inv-1}
x_{i}=q_{i}(\theta_{l},\theta_{2},\ldots,\theta_{n}).
\end{equation}
Для сокращения в записи введем обозначения:
$
\overrightarrow{\theta}=
\left\{\theta_{l},\theta_{2},\ldots,\theta_{n} \right\}$,
$\overrightarrow{q}=\left\{q_{l},q_{2},\ldots,q_{n} \right\}
\in {\mathbb R}^{n},
$
$
\overrightarrow{\rho}(t;{\bf x})=\left\{\rho_{1}(t;{\bf
x}),\rho_{2}(t;{\bf x}),\ldots, \rho_{3}(t;{\bf x})   \right\}
\in {\mathbb R}^{n}.
$
Так как в силу условия \eqref{inv-1} для некоторого $s\geq1$ и любой другой функции
$
\chi(t;{\bf x})=\rho_{n+s}(t;{\bf x})\,\rho_{n+1}^{-1}(t;{\bf x})
$
в момент времени $t=0$ возможно представление
$$
\chi(0;{\bf y})=\rho_{n+s}(0;{\bf y})\,\rho_{n+1}^{-1}(0;{\bf y})=\overline{\psi} \Bigl(
\overrightarrow{q}(\overrightarrow{\theta})\Bigr),$$
то, следовательно, для любой функции $\rho_{n+s}(t;{\bf x})$,
являющейся ядром интегрального инварианта системы \eqref{Ayd01.1} в силу
утверждения леммы \ref{l1}, получаем, с учетом построения функций $\theta_{l}(t;{\bf x})$:
$$
\chi(t;{\bf x}(t;{\bf  y}))=\rho_{n+s}(t;{\bf x}(t;{\bf
y}))\,\rho_{n+1}^{-1}(t;{\bf x}(t;{\bf  y}))=
$$
$$
=\rho_{n+s}(0;{\bf y})\,\rho_{n+1}^{-1}(0;{\bf y})=\chi(0;{\bf y})=
\overline{\psi} \biggl(
\overrightarrow{q}\Bigl(\rho_{n+1}^{-1}(t;{\bf x}(t;{\bf  y}))\cdot\overrightarrow{\rho}(t;{\bf x}(t;{\bf
y})) \Bigr)\biggr).
$$
Из этого соотношения и следует, что для любого $ t\geq 0 $ и всех $s\geq1$
$$
\rho_{n+s}(t;{\bf x}(t;{\bf  y}))=\rho_{n+1}(t;{\bf x}(t;{\bf
y}))\, \overline{\psi}  \biggl(
\overrightarrow{q}\Bigl(\rho_{n+1}^{-1}(t;{\bf x}(t;{\bf  y}))\cdot\overrightarrow{\rho}(t;{\bf x}(t;{\bf
y})) \Bigr)\biggr).
$$
Т. е. полная совокупность ядер состоит из $(n+1)$-й линейно независимых функций.
\end{proof}

\begin{cor}\label{c1}
Полная совокупность линейно независимых стохастических первых интегралов системы \eqref{Ayd01.1} состоит из $n$ функций.
\end{cor}
\begin{proof}
Функция
$
\widetilde{u}_{l}(t;{\bf x})=\rho_{l}(t;{\bf x})\rho_{n+1}^{-1}(t;{\bf x})
$
обладает свойствами, описанными в определении \ref{Ayddf1} и является стохастическим первым интегралом системы ДУ \eqref{Ayd01.1}. Таких функций получено ровно $n$.
\end{proof}

\begin{remark}
Пусть случайный процесс ${\bf x}(t)$ является
решением обобщенного уравнения Ито, которое представим в виде:
\begin{equation*}
    {d}_{t}{\bf x}(t)=a(t;{\bf x}(t))dt + b(t;{\bf
    x}(t))d{\bf w}(t)+ dP(t,\Delta\gamma)=\widetilde{d}_{t}{\bf x}(t)+dP(t,\Delta\gamma).
\end{equation*}
Тогда обобщенную формулу Ито \eqref{Ayd2} можно записать в виде:
\begin{equation}\label{AydIt2}
       d_{t}f(t;{\bf x}(t))= \widetilde{d}_{t}f(t;{\bf
       x}(t))+\widetilde{d}_{t}P(t,\Delta\gamma),
\end{equation}
где $\widetilde{d}_{t}f(t;{\bf x}(t))$ -- дифференциал Ито для
части $\widetilde{d}_{t}{\bf x}(t)$, $\widetilde{d}_{t}P(t,\Delta\gamma)$ --
дифференциал пуассоновской добавки $dP(t,\Delta\gamma)$.
\end{remark}

Построим   уравнение для $u(t;{\bf x})$, воспользовавшись следствием \ref{c1} и соотношением:
\begin{equation}\label{Ayd2.11}
\ln u_{s}(t;{\bf x})=\ln \rho_{s}(t;{\bf x})-\ln\rho_{l}(t;{\bf
x}).
\end{equation}

Продифференцируем $\ln \rho(t;{\bf x})$ с использованием
\eqref{Aydii} и обобщенной формулы Ито и (\ref{AydIt2}).
\begin{equation}\label{Aydln}
\begin{array}{c}
  \displaystyle d_{t}\ln\rho(t;{\bf x})= \frac{1}{\rho(t;{\bf
x})}\widetilde{d}_{t}\rho(t;{\bf x})-\frac{1}{2\rho^{2}(t;{\bf
x})}\left(-\displaystyle\frac{\partial(\rho(t;{\bf
x})b_{i\,k}(t;{\bf
x})) }{\partial x_{i}}\right)^{2}dt +  \\
  +\displaystyle\int\limits_{R(\gamma)}\Bigl[\Bigr.\ln\left\{\right.
\rho_{s}\left( t;{\bf x}-g(t;{\bf
  x}(t;{\bf y});\gamma);\gamma\right)
  \cdot D\left( {\bf x}^{-1}(t;{\bf x};\gamma)  \right)
  \left.\right\}
  -\ln\rho_{s}(t;{\bf x})\Bigl.\Bigr]\nu(dt;d\gamma).
\end{array}
\end{equation}
Рассмотрим сумму первых двух слагаемых, для компактности записи
опуская аргумент у функций:
\begin{equation}\label{Aydlnp}
\begin{array}{c}
  S_{1}=\displaystyle \frac{1}{\rho}\, \widetilde{d}_{t}\rho-\frac{1}{2\rho^{2}}\,
  \left(-\displaystyle\frac{\partial(\rho b_{i\,k})}{\partial x_{i}}\right)^{2} =\\
  = \displaystyle\frac{1}{\rho} \biggl[\Bigl( -\frac{\partial (\rho a_{i})}{\partial x_{i}}+
  \frac{1}{2}\frac{\partial ^{\,2}(\rho b_{i\,k}b_{j\,k})}{\partial x_{i}\partial x_{j}}
  \Bigr)dt + \Bigl( -\displaystyle\frac{\partial(\rho b_{i\,k}) }{\partial
  x_{i}}\Bigr)dw_{k}(t)\biggr]
  -\\
  -\displaystyle\frac{1}{2\rho^{2}}\,\left(-\displaystyle\frac{\partial(\rho b_{i\,k})}{\partial
 x_{i}}\right)^{2}dt=\\
 = \displaystyle\biggl[\biggr.  - \frac{\partial  a_{i}}{\partial x_{i}}-a_{i}\frac{\partial \ln\rho}{\partial x_{i}} +
 \frac{1}{2\rho}\frac{\partial}{\partial x_{i} }\left(\rho\,\frac{\partial b_{i\,k}b_{j\,k}}{\partial
  x_{j}}+ b_{i\,k}b_{j\,k}\,\frac{\partial \rho}{\partial x_{i}} \right)-
  \\
  -\displaystyle\frac{1}{2\rho^{2}}\,\left(\displaystyle\rho\,\frac{\partial b_{i\,k}}{\partial
 x_{i}}+b_{i\,k}\, \frac{\partial \rho}{\partial
 x_{i}}\right)^{2}\biggl.\biggr]dt
 -\left( b_{i\,k}\displaystyle\frac{\partial\ln\rho }{\partial
 x_{i}}+\frac{\partial b_{i\,k}}{\partial
 x_{i}}\right)dw_{k}(t)=\\
 =S_{2}dt-\left( b_{i\,k}\displaystyle\frac{\partial\ln\rho }{\partial
 x_{i}}+\frac{\partial b_{i\,k}}{\partial
 x_{i}}\right)dw_{k}(t).
\end{array}
\end{equation}
Преобразуем  часть $S_{2}$:
\begin{equation}\label{Aydlnp1}
\begin{array}{c}
S_{2}= - \displaystyle \frac{\partial  a_{i}}{\partial
x_{i}}-a_{i}\frac{\partial \ln\rho}{\partial x_{i}} +
\frac{\partial b_{i\,k}b_{j\,k}}{\partial
  x_{j}}\frac{\partial \ln\rho}{\partial x_{i}}  +
  \frac{1}{2}\,\frac{\partial^{\,2}( b_{i\,k}b_{j\,k})}{\partial x_{i}\partial x_{j}}
  +\frac{1}{2}b_{i\,k}b_{j\,k}\,\frac{1}{\rho}\,\frac{\partial^{\,2} \rho}{\partial x_{i}\partial x_{j}}
  -\\
  -
  \displaystyle\frac{1}{2}\,\left(\displaystyle\frac{\partial b_{i\,k}}{\partial
 x_{i}}+b_{i\,k}\, \frac{\partial \ln\rho}{\partial
 x_{i}}\right)\left(\displaystyle\frac{\partial b_{j\,k}}{\partial
 x_{j}}+b_{j\,k}\, \frac{\partial \ln\rho}{\partial
 x_{j}}\right)=\\
 =- \displaystyle \frac{\partial  a_{i}}{\partial
x_{i}}-a_{i}\frac{\partial \ln\rho}{\partial x_{i}} +
\frac{\partial b_{i\,k}b_{j\,k}}{\partial
  x_{j}}\frac{\partial \ln\rho}{\partial x_{i}}  +
  \frac{1}{2}\,\frac{\partial^{\,2}( b_{i\,k}b_{j\,k})}{\partial x_{i}\partial
  x_{j}}+\frac{1}{2}\,b_{i\,k}b_{j\,k}\frac{\partial^{\,2} \ln\rho}{\partial x_{i}\partial
  x_{j}}+\\
  +\displaystyle\frac{1}{2}\,b_{i\,k}b_{j\,k}\frac{\partial\ln\rho}{\partial x_{i}}\,\frac{\partial\ln\rho}{\partial
  x_{j}}
  -\displaystyle\frac{1}{2}\left[\frac{\partial b_{i\,k}}{\partial x_{i}}\,
    \frac{\partial b_{j\,k}}{\partial x_{j}} +2\, b_{i\,k}\frac{\partial b_{j\,k}}{\partial
    x_{j}}\,\frac{\partial \ln\rho}{\partial x_{i}} +
    b_{i\,k}b_{j\,k}\,\frac{\partial\ln\rho}{\partial x_{i}}\,\frac{\partial\ln\rho}{\partial
  x_{j}}  \right]=\\
  =- \displaystyle \frac{\partial  a_{i}}{\partial
x_{i}}-a_{i}\frac{\partial \ln\rho}{\partial x_{i}} +
\frac{\partial b_{i\,k}b_{j\,k}}{\partial
  x_{j}}\,\frac{\partial \ln\rho}{\partial x_{i}}  +
  \frac{1}{2}\,\frac{\partial^{\,2}( b_{i\,k}b_{j\,k})}{\partial x_{i}\partial
  x_{j}}+\frac{1}{2}\,b_{i\,k}b_{j\,k}\frac{\partial^{\,2} \ln\rho}{\partial x_{i}\partial
  x_{j}}-\\
    -\displaystyle\frac{1}{2}  \frac{\partial b_{i\,k}}{\partial x_{i}}\,
    \frac{\partial b_{j\,k}}{\partial x_{j}} - b_{i\,k}\frac{\partial b_{j\,k}}{\partial
    x_{j}}\,\frac{\partial \ln\rho}{\partial x_{i}}.
   \end{array}
\end{equation}
Следовательно, подставляя (\ref{Aydlnp1}) в (\ref{Aydlnp}), затем
в (\ref{Aydln}), получаем следующее уравнения для  $\ln
\rho(t;{\bf x})$:
\begin{equation}\label{Ayd2.11a}
\begin{array}{c}
  d_{t}\ln\rho(t;{\bf x})=\biggl[\biggr.- \displaystyle
  \frac{\partial a_{i}(t;{\bf x})}{\partial x_{i}}-
  a_{i}(t;{\bf x})\frac{\partial \ln\rho(t;{\bf x})}{\partial x_{i}}
  -\displaystyle \frac{1}{2}\,\frac{\partial b_{i\,k}(t;{\bf x})}
  {\partial x_{i}}\,\frac{\partial b_{j\,k}(t;{\bf x})}
  {\partial x_{j}}+\\
  +
  \displaystyle \frac{1}{2}b_{i\,k}(t;{\bf x})b_{j\,k}(t;{\bf x})
  \frac{\partial^{2} \ln\rho(t;{\bf x})}{\partial x_{i}\partial x_{j}}
  +\displaystyle\frac{\partial \left(b_{i\,k}(t;{\bf x})b_{j\,k}(t;{\bf x})\right)}
  {\partial x_{i}}\,
  \frac{\partial \ln\rho(t;{\bf x})}{\partial  x_{j}}+\\
  +\displaystyle\frac{1}{2}\frac{\partial^{2}\left( b_{i\,k}(t;{\bf x})b_{j\,k}(t;{\bf x})\right)}
  {\partial x_{i}\partial x_{j}}
  - b_{j\,k}(t;{\bf x})\displaystyle\frac{\partial b_{i\,k}(t;{\bf x})}
  {\partial x_{i}}\,\frac{\partial \ln\rho(t;{\bf x})}{\partial x_{j}}\biggl.\biggr] dt +\\
+\displaystyle\int\limits_{R(\gamma)}\Bigl[\Bigr.\ln\left\{\right.
\rho\left( t;{\bf x}-g(t;{\bf
  x}(t;{\bf y});\gamma);\gamma\right)
  \cdot D\left( {\bf x}^{-1}(t;{\bf x};\gamma)  \right)
  \left.\right\}-\\
  -\ln\rho(t;{\bf x})\Bigl.\Bigr]\nu(dt;d\gamma)
  -\biggl[\displaystyle \frac{\partial b_{i\,k}(t;{\bf x})}{\partial x_{i}}+
b_{j\,k}(t;{\bf x})\frac{\partial \ln\rho(t;{\bf x})}{\partial
x_{j}}
   \biggr]dw_{k}(t).
\end{array}
\end{equation}

Опираясь на полученное уравнение (\ref{Ayd2.11a}), составим
уравнения для $\rho_{s}(t;{\bf x})$ и $\rho_{l}(t;{\bf x})$ и с
учетом (\ref{Ayd2.11}) приходим к
выводу, что стохастический первый интеграл $u(t;{\bf x};\omega)$
обобщенного уравнения Ито является решением уравнения:
\begin{equation}\label{Ayd2.12}
\begin{array}{c}
  d_{t}u(t;{\bf x})=\biggl[\biggr. -  a_{i}(t;{\bf x})\displaystyle\frac{\partial u(t;{\bf x})}{\partial
x_{i}} +\displaystyle \frac{1}{2}\,b_{i\,k}(t;{\bf
x})b_{j\,k}(t;{\bf x})
  \frac{\partial^{2} u(t;{\bf x})}{\partial x_{i} \partial x_{j}} -\\
-  b_{i\,k}(t;{\bf x})\displaystyle\frac{\partial}{\partial
x_{i}}\left(b_{j\,k}(t;{\bf x})\frac{\partial u(t;{\bf
x})}{\partial x_{j}}\right)\biggl.\biggr]dt
  - b_{i\,k}(t;{\bf x})\displaystyle\frac{\partial u(t;{\bf x})}{\partial
x_{i}}\,dw_{k}(t)+ \\
+ \displaystyle\int\limits_{R(\gamma)}\Bigl[u\Bigl(t; {\bf
x}-g(t;{\bf
  x}^{-1}(t;{\bf x};\gamma));\gamma\Bigr)
 -
  u(t;{\bf x})\Bigr]\nu(dt;d\gamma).
\end{array}
\end{equation}

В работе {\rm{\cite{D_78}}} было введено  понятие первого
интеграла для системы стохастических дифференциальных уравнений
Ито, как неслучайных функций на произвольных случайных
реализациях возмущений:
$$
u\Bigl(0;{\bf x}(0)\Bigr)=u\Bigl(t;{\bf x}(t; {\bf x}(0))\Bigr).
$$
Однако можно говорить и о стохастическом первом интеграле для системы стохастических дифференциальных уравнений Ито. Следует отметить, что
при отсутствии пуассоновских возмущений также можно говорить о {\it стохастическом} (а не только детерминированном) первом интеграле.

Для случая только винеровских возмущений (классическое уравнение
Ито) в {\rm{\cite{D_78}}} исследование свойств первого интеграла как детерминированной функции было связано с установлением независимости такой функции
от реализаций ${\bf w}(t)$.
В рассматриваемом нами в рамках определения \ref{Ayddf1} случае функция  $u\Bigl(t;{\bf x}\Bigr)$, для которой
построено уравнение (\ref{Ayd2.12}), зависит также и от
пуассоновских возмущений, т. е. является случайной.

При добавлении
пуассоновских возмущений (обобщенное уравнение Ито) такое
требование, как следует из уравнения (\ref{Ayd2.12}), приводит к
следующим условиям.
\begin{theorem}\label{th-usl}
Случайная функция $u\left(t;{\bf x}(t)\right) \in \mathcal{C}_{t,x}^{1,2}$, определенная на том же вероятностном пространстве, что и случайный процесс ${\bf x}(t)$, являющийся решением системы обобщенных стохастических дифференциальных уравнений Ито \eqref{Ayd01.1}, есть первый интеграл системы \eqref{Ayd01.1} тогда и только тогда, когда функция  $u\left(t;{\bf x}(t)\right)$ удовлетворяет условиям $\left.\mathcal{L}\right)$\label{uslL1}:
\begin{enumerate}
    \item $b_{i\,k}(t;{\bf x})\displaystyle\frac{\partial u(t;{\bf x})}{\partial
x_{i}}=0$, для всех $k=\overline{1,m}$ (компенсация винеровского
возмущения);
    \item $\displaystyle\frac{\partial u(t;{\bf x})}{\partial
t}+\displaystyle\frac{\partial u(t;{\bf x})}{\partial
x_{i}}\Bigl[a_{i}(t;{\bf x})-
\displaystyle\frac{1}{2}\,b_{j\,k}(t;{\bf x})\frac{\partial
b_{i\,k}(t;{\bf x})}{\partial x_{j}}\Bigr]=0$ (независимость от
времени);
     \item $u(t;{\bf x})-u\Bigl(t;{\bf x}+g(t;{\bf
    x};\gamma)\Bigr)=0$ для любых $\gamma\in R(\gamma)$
         во всей области
  определения процесса (компенсация пуассоновских скачков).
\end{enumerate}
\end{theorem}

\begin{proof}
Рассмотрим уравнение \eqref{Ayd2.12}, решением которого является функция
$u\left(t;{\bf x}(t)\right)$. Перенесем все в одну сторону от знака равенства. Учитывая, что $d_{t}u(t;{\bf x})=\dfrac{\partial u(t;{\bf x})}{\partial t}\,dt$, получаем для любого $t$:
\begin{equation*}
\begin{array}{c}
  \biggl[\biggr.\dfrac{\partial u(t;{\bf x})}{\partial t}+   a_{i}(t;{\bf x})\displaystyle\frac{\partial u(t;{\bf x})}{\partial
x_{i}} -\displaystyle \frac{1}{2}\,b_{i\,k}(t;{\bf
x})b_{j\,k}(t;{\bf x})
  \frac{\partial^{2} u(t;{\bf x})}{\partial x_{i} \partial x_{j}} +\\
+   b_{i\,k}(t;{\bf x})\displaystyle\frac{\partial}{\partial
x_{i}}\left(b_{j\,k}(t;{\bf x})\frac{\partial u(t;{\bf
x})}{\partial x_{j}}\right)\biggl.\biggr]dt
  + b_{i\,k}(t;{\bf x})\displaystyle\frac{\partial u(t;{\bf x})}{\partial
x_{i}}\,dw_{k}(t)- \\
- \displaystyle\int\limits_{R(\gamma)}\Bigl[u\Bigl(t; {\bf
x}-g(t;{\bf
  x}^{-1}(t;{\bf x};\gamma));\gamma\Bigr)
 +
  u(t;{\bf x})\Bigr]\nu(dt;d\gamma)=0.
\end{array}
\end{equation*}
Следовательно, множители при $dt$, при $dw(t)$ и при $\nu(dt;d\gamma)$ должны быть нулевыми. "<Винеровское"> слагаемое равно нулю:
\begin{equation}\label{usl-win}
b_{i\,k}(t;{\bf x})\displaystyle\frac{\partial u(t;{\bf x})}{\partial
x_{i}}=0 \ \ \ \ \textrm{для всех}\ \ k=\overline{1,m}.
\end{equation}
Для "<пуассоновской"> части получаем:
\begin{equation}\label{usl-puass}
u\Bigl(t; {\bf x}-g(t;{\bf x}^{-1}(t;{\bf x};\gamma));\gamma\Bigr)-
u(t;{\bf x})=0.
\end{equation}
Преобразуем равенство \eqref{usl-puass}. Перейдем к переменным
$${\bf y}={\bf x}-g(t;{\bf x}^{-1}(t;{\bf x};\gamma);\gamma).$$
Учитывая, что ${\bf x}^{-1}(t;{\bf x};\gamma)$ -- обозначение
обратной функции для $f({\bf x})={\bf x}+g(t;{\bf x};\gamma)$ (см.
(\ref{Aydx1y})), убеждаемся, что это условие эквивалентно
следующему:
\begin{equation}\label{usl-puass-1}
u(t;{\bf x})-u\Bigl(t;{\bf x}+g(t;{\bf
    x};\gamma)\Bigr)=0 \ \ \ \textrm{для любых} \ \ \ \gamma\in R(\gamma).
\end{equation}
Далее используем правило дифференцирование произведения
 и условие \eqref{usl-win}, получаем:
\begin{equation*}
\begin{array}{c}
  \dfrac{\partial u(t;{\bf x})}{\partial t}+
  a_{i}(t;{\bf x})\dfrac{\partial u(t;{\bf x})}{\partial
x_{i}} + \dfrac{1}{2}\,b_{i\,k}(t;{\bf
x})b_{j\,k}(t;{\bf x})
  \dfrac{\partial^{2} u(t;{\bf x})}{\partial x_{i} \partial x_{j}}\, -\\
-   b_{i\,k}(t;{\bf x})\dfrac{\partial}{\partial
x_{i}}\left(b_{j\,k}(t;{\bf x})\dfrac{\partial u(t;{\bf
x})}{\partial x_{j}}\right)= \\
  =  \dfrac{\partial u(t;{\bf x})}{\partial t}+
  a_{i}(t;{\bf x})\cfrac{\partial u(t;{\bf x})}{\partial
x_{i}}\, +\\
+ \dfrac{1}{2}\,b_{i\,k}(t;{\bf
x})\left[\dfrac{\partial}{\partial x_{i}}\Bigl(b_{j\,k}(t;{\bf x})
\dfrac{\partial u(t;{\bf x})}{\partial x_{j}}\Bigr)-
\dfrac{\partial u(t;{\bf x})}{\partial x_{j}}
\dfrac{\partial b_{jk}(t;{\bf x})}{\partial x_{i}}\right]=\\
= \dfrac{\partial u(t;{\bf x})}{\partial t}+
  a_{i}(t;{\bf x})\cfrac{\partial u(t;{\bf x})}{\partial
x_{i}}\, - \dfrac{1}{2}\,b_{i\,k}(t;{\bf
x})
\dfrac{\partial u(t;{\bf x})}{\partial x_{j}}
\dfrac{\partial b_{jk}(t;{\bf x})}{\partial x_{i}}.
\end{array}
\end{equation*}
Следовательно,
$$\displaystyle\frac{\partial u(t;{\bf x})}{\partial
t}+\displaystyle\frac{\partial u(t;{\bf x})}{\partial
x_{i}}\Bigl[a_{i}(t;{\bf x})-
\displaystyle\frac{1}{2}\,b_{j\,k}(t;{\bf x})\frac{\partial
b_{i\,k}(t;{\bf x})}{\partial x_{j}}\Bigr]=0.$$
Таким образом, получили все условия $\left.\mathcal{L}\right)$.
\end{proof}

\begin{remark}
В случае, когда рассматриваем конкретную реализацию, т. е.
параметр $\omega$ в дальнейшем не влияет, неслучайную функцию
$u(t;{\bf x})$ можно считать детерминированным первым интегралом
стохастической системы.
\end{remark}

\begin{remark}В {\rm\cite[c.\,24]{D_02}} было введено понятие стохастического первого
интеграла для центрированной пуассоновской меры и полученные
условия для его существования учитывают необходимость задания
плотности интенсивности пуассоновского распределения в отличие от
предложенного в данной статье. Таким образом, безразлично, каков
вероятностный закон имеют интенсивности пуассоновских скачков. Это
обстоятельство является очень важным для дальнейших применений, в
частности, для построения программных управлений {\rm\cite{11_KchUpr}}.
\end{remark}

\begin{remark}В {\rm\cite{Oks_07}} была предложена формула, являющаяся обобщением  формулы Ито-Вентцеля для СДУ Ито, содержащего только пуассоновскую составляющую с центрированной пуассоновской мерой.
\end{remark}

Предложенное обобщение формулы Ито-Вентцеля  и понятия
стохастического первого интеграла {\rm{\cite{D_02}}} позволяет,
как отмечено в {\rm{\cite{07_ChOboz}}}, строить программные
управления стохастических динамических систем, подверженных
случайным возмущениям, вызванным винеровскими возмущениями и
пуассоновскими скачками
{\rm{\cite{11_KchUpr}}}.

\end{document}